\documentclass{amsart}
\usepackage{amsmath,amssymb}
\usepackage{url}

\newcommand{\nats}{\mathbb{N}}

\newtheorem{lemma}{Lemma}
\newtheorem{theorem}{Theorem}
\newtheorem{corollary}{Corollary}
\newtheorem{proposition}{Proposition}
\newcommand{\MFN}{{\sf MFn}}

\begin{document}
\author[Thierry Coquand and Bas Spitters]
         {THIERRY COQUAND\\
          Computing Science Department at G\"oteborg University\\
          BAS SPITTERS\\
          Department of Mathematics and Computer Science, Eindhoven University
of Technology}

\title{Constructive Gelfand duality for C*-algebras}\maketitle

\begin{abstract}
  We present a constructive proof of Gelfand duality for C*-algebras by
  reducing the problem to Gelfand duality for real C*-algebras.
\end{abstract}

\section{Introduction}

Classical Gelfand duality states that category of commutative C*-algebras and
the category of compact Hausdorff spaces are equivalent. The proof relies on
the axiom of choice in an essential way.
In a sequence of papers starting in
a 1980 pre-print and culminating in the references
\cite{banaschewskimulvey00a,banaschewskimulvey00b,BMStone,banaschewskimulvey06},
Banaschewski and Mulvey explore
a constructive version of the Gelfand duality theorem which can be applied
internally in
a topos. In this context, the category of compact Hausdorff spaces
is replaced by the category of compact completely regular locales. A locale is
is a pointfree topology: a lattice theoretic presentation of the open sets of a
topological space. In the presence of the axiom of choice, the category of
compact completely regular locales and the category of compact Hausdorff spaces
are equivalent. The axiom of choice is (only) used to
construct the points in the topological spaces. In topos theory, the axiom of
choice is not
generally present~\cite{Mulvey:geometry}. In this light,
Banachewski and Mulvey generalized Gelfand
duality to Grothendieck toposes by rephrased it as the equivalence of the
category of commutative C*-algebras and the category of compact completely
regular locales.
When the axiom of choice is present the spatial version is a simple corollary.

The treatment by Banachewski and Mulvey is not quite constructive: it relies on
Barr's Theorem. Barr's theorem states: If a geometric statement is
deducible from a geometric theory using classical logic and the axiom of
choice, then it is also deducible from it constructively; see~\cite{Wraith} for
a discussion of the importance of this theorem in constructive algebra. The
proof of Barr's theorem itself, however, is highly non-constructive. Even if we
are willing to grant this, Barr's theorem depends on the topos being a
Grothendieck
topos.

We give a fully
constructive treatment of Gelfand duality. An alternative constructive proof of
Gelfand duality is
announced
in~\cite{banaschewskimulvey06} and~\cite{Mulvey:geometry}.
Our proof uses a concrete
presentation of the Gelfand spectrum as a lattice. Such constructive proofs are
sometimes more direct~\cite{coquandspitters05} than proofs via an encoding of
topology in metric spaces, as is common in Bishop's constructive
mathematics~\cite{Bishop67}. Moreover, this construction of the lattice
presenting the spectrum as a locale is technically advantagous, as it is
preserved under inverse images of geometric morphisms. As such it has been
applied in~\cite{qtopos}.

The article is organized as follows. We start by a constructive reduction of
Gelfand
duality from
the complex case to the real case (Section~\ref{reduction}). A constructive
presentation of Gelfand duality in the
real case has been given in~\cite{coquand05}. In order to apply these results we
prove that the set of self-adjoint elements of a C*-algebra is a real
C*-algebra (Section~\ref{realcstar}). We put all the pieces together in
Section~\ref{gelfand}. Finally, Section~\ref{applications} ends with short
direct proofs of results which were obtained using Barr's theorem
in~\cite{banaschewskimulvey06}.

\section{Preliminaries}

 We recall here the definition of  a commutative $C^*$-algebra $A$ in a
topos following~\cite{banaschewskimulvey06}.
 When working in an intuitionistic framework, we cannot assume
in general the (semi)norm of an element to be a Dedekind real, but instead it
may simply
be a {\em non negative upper real}.
We define a non negative upper
real to be a inhabited open upward closed set of positive rational numbers. We
can define the
addition and multiplication of non negative upper reals: $U_1+U_2$ is the set of rationals
$r_1+r_2,~r_1\in U_1,~r_2\in U_2$ and $U_1U_2$ is the set of rationals
$r_1r_2,~r_1\in U_1,~r_2\in U_2$. We define also $U_1\leqslant U_2$ to mean that
$U_2$ is a subset
of $U_1$. Finally we may identify the non negative rational $q$ with the set of
rationals
$r$ such that $r>q$. The norm $\| a\|$ of $a$ in $A$ is then an upper real.
The notation of \cite{banaschewskimulvey06} is $a\in N(q)$ for $\|a\|<q$.
The conditions
for the relation $a\in N(q)$, introduced in~\cite{banaschewskimulvey06} can then
be written as the
usual conditions on the seminorm
$$
\|0\| = 0,~~~~\|1\| = 1,~~~~\|a^*\| = \|a\|,~~~~~\|ab\|\leqslant \|a\| \|b\|
$$
$$
\|ra\| = |r| \|a\|,~~~~~~\|a+b\| \leqslant \|a\|+\|b\|,~~~~~\|aa^*\| = \|a\|^2
$$
As in \cite{banaschewskimulvey06}, we assume finally $A$ to be {\em complete}:
any Cauchy
approximation on $A$ has a unique limit in $A$. (As a consequence, $a=0$ iff
$\|a\|=0$).

We will use the letters $a,b,x,y$ to range over elements of the C*-algebra and
the letters $q,r,s,t$ to range over the rationals.

\section{\label{reduction}Reduction to the real case}
Let $A$ be a C*-algebra and let $B = A_{sa}$ be the set of {\em self-adjoint}
elements, i.e.\ elements $a$ such that $a^* = a$.
The algebra $B$ is then a commutative Banach algebra over the rationals.
For $a$ in $B$, we have $\|a^2\|=\|a\|^2$, since $a=a^\ast$.

\begin{proposition}\label{selfadjoint}
 For $a,b$ in $B$ we have $\|a^2 \| \leqslant \|a^2 + b^2 \|$.
\end{proposition}

\begin{proof}
  We write $a^2 + b^2 = (a + b i)
  (a - b i) =: c c^{\ast}$. So $\|a^2 + b^2 \|=\|c c^{\ast} \|=\|c\|^2$.
  Finally, $2 a = c + c^{\ast}$, so $\|a\|= \frac{1}{2} \|c + c^{\ast} \|
  \leqslant \frac{1}{2} (\|c\|+\|c^*\|) =\|c\|$ and therefore $\|a^2 \|=\|a\|^2
  \leqslant \|c\|^2 =\|a^2 + b^2 \|$.
\end{proof}

\section{\label{realcstar}Real Banach algebras}

 In this section, we consider a complete commutative Banach algebra
$B$ over the rationals such that $\|a^2\|=\|a\|^2$ and $\|a^2 \| \leqslant \|a^2
+ b^2 \|$.
By Proposition~\ref{selfadjoint}, this will be the case if we take for $B$ the
self-adjoint part of a commutative $C^*$-algebra.

\begin{lemma}
  \label{lem:square}If $\|1 - x\| \leqslant 1$. Then $x$ is a square.
\end{lemma}

\begin{proof}We give an explicit proof that the Taylor series for
$\sqrt{1-(1-x)}$ converges.
  We define two sequences: $y_n$ in $B$ and $r_n$ in $\mathbb{Q}$. We take $y_0
= 0,~r_0=0$ and
$y_{n + 1} = \frac{1}{2} (1 - x + y^2_n)$ and $r_{n + 1} =  \frac{1}{2} (1 + r_n^2)$.

  For all $n$, $\|y_n \| \leqslant r_n$ by induction. Since we have
  \[y_{n + 1} - y_n = \frac{1}{2} (y_n + y_{n - 1}) (y_n - y_{n - 1})\]
  we get $\|y_{n + 1} - y_n \| \leqslant r_{n + 1} - r_n$ by induction.
  Consequently, \[\|(1 - y_n)^2 - x\|= 2\|y_{n + 1} - y_n \| \leqslant 2 (r_{n
+ 1} - r_n) \rightarrow 0\]
because we have $r_n\rightarrow 1$ in a constructive way \cite{coquand05}.
\end{proof}

\begin{proposition}\label{sum:square}
  A sum of squares is a square.
\end{proposition}

\begin{proof}
  As in \cite{Kelley}. We claim that $\|x\|, \|1 - x\| \leqslant 1$ iff $x$
and $1 - x$ are squares.

The implication from left to right is Lemma~\ref{lem:square}. For the reverse
implication suppose that $x = u^2$ and $1 - x = v^2$, then $1 = u^2 + v^2$, so
  $\|u\|^2, \|v\|^2 \leqslant 1$.

  For the proof of the Proposition let $x, y$ be squares. We can
assume $\|x\|, \|y\| \leqslant 1$. Then $1 - x$ and $1 - y$ are squares and so
$\|1 - x\|, \|1 - y\| \leqslant 1$.
  Since \[
\|1 - \frac{(x + y)}{2} \| \leqslant \frac{1}{2} (\|1 - x\|+\|1 -
  y\|) \leqslant 1,
        \]
$(x + y) / 2$ is a square and so is $x + y$.
\end{proof}

 Let $P$ be the set of all squares. Then $P$ is a {\em cone}: it contains the squares
and is closed under multiplication and addition. The cone $P$ defines an
ordering on the algebra $B$. As in \cite{coquand05} we define $r\ll a$ to mean
$a-s\in P$
for some $s>r$.
By Lemma~\ref{lem:square} we have $r-a$ in $P$ if $\|a\|\leqslant r$ and hence
$B$ has the multiplicative unit 1 as a strong unit for this ordering.
Consequently, all the results of the
first part of
\cite{coquand05} are available.

 We define $\MFN (B)$ to be the locale generated by symbols $D(a)$, $a \in B$,
and relations
\begin{enumerate}
\item $D(1) = 1$
\item $D(-a^2) = 0$
\item $D(a+b)\leqslant D(a)\vee D(b)$
\item $D(a)\wedge D(-a) = 0$
\item $D(ab) = (D(a)\wedge D(b))\vee (D(-a)\wedge D(-b))$
\item $D(a) = \bigvee_{r>0} D(a-r)$
\end{enumerate}

The points of this locale are the Multiplicative Functionals. A symbol $D(a)$
intuitively represents the open set $\{\phi \mid \phi(a)>0\}$.

\begin{lemma}\label{key}
If $0\ll ac$ and $0\leqslant c$ then $0\ll a$.
\end{lemma}

\begin{proof}
See \cite{Krivine} Th\'eor\`eme 12. We give a sketch of the argument.
Since the ring is Archimedean, we have $N$ in $\nats$ such that
$-N\leqslant a\leqslant N$.
Since $0\leqslant c$ and $1\leqslant ac$ we have $1\leqslant Nc$ and thus
$\frac1N\leqslant c$. There exists $L$ in $\nats$ such that $c\leqslant L$ and
we get $\frac1N\leqslant c\leqslant L$. If we write $b = 1-\frac{c}L$, we have
$0\leqslant b\leqslant 1 -\frac1{NL}$
and $\frac1L\leqslant a(1-b)$. By multiplying
by $1+\dots+b^{n-1}$ we get $\frac1L\leqslant a(1-b^n)$ and so
$\frac1L+ab^n\leqslant a$.
For $n$ big enough we have $b^n\leqslant \frac1{2NL}$; hence
$\frac1{2L}\leqslant a$.
\end{proof}

 One of the main results of \cite{coquand05} is a constructive proof of the
following result.

\begin{proposition}\label{main}
We have $D(a) =1$ in $\MFN (B)$ iff $0\ll a$ in $B$.
\end{proposition}
\begin{proof}
The proof which we sketch here is a combination of Lemma \ref{key} and a
cut-elimination argument
\cite{coquand:entail,DynMethod}, which is an important technique in proof
theory.

First we derive some simple consequences of the axioms (1-5).
\begin{itemize}
 \item[-] If $a\leqslant b$, that is, $b-a\in P$, then $D(a)\leqslant D(b)$:\\
$a=b-x^2$, so
$D(a)\leqslant  D(b)\vee D(-x^2)$, which is equal to $D(b)\vee 0=D(b)$.
 \item[-] For all $n$, $D(\frac1n)=1$, from (1) and (3).
\end{itemize}
It follows that we have $D(s)=1$ if $s>0$ and that $D(a)=1$ if $0\ll a$.
This is the implication from right to left.

We now consider the converse direction.

First we notice that $D(a)=1$ follows from (1-5) iff it follows from (1-6). For
this we define an interpretation of the theory (1-6) into (1-5) by
reinterpreting the symbol $D(a)$ as $\bigvee_{r>0}D(a-r)$;
see~\cite{banaschewskimulvey00a, coquand05}.

Next, we characterise the distributive lattice generated by (1-5). We have \[
D(a_1)\wedge \ldots \wedge D(a_n)\leqslant  D(b_1)\vee \ldots \vee D(b_m)
                                                                           \]
iff we have a relation $m+p=0$, where $m$ belongs to the multiplicative monoid
generated by $a_1,\ldots,a_n$ and $p$ belongs to the $P$-cone generated by
$-b_1,\ldots,-b_m$. A \emph{$P$-cone} is a subset which contains $P$ and
is closed under addition and multiplication. For the proof see~\cite{DynMethod,
coquand05}.

It follows that if $D(a)=1$ in (1-5), then we have a relation $m+p=0$, where
$m=1$ and $p$ belongs to the $P$-cone generated by $-a$. Hence, there are $b,c$
in $P$ such that $1+b+c(-a)=0$, that is $ca=1+b$. Consequently, $0\ll a$ by
Lemma~\ref{key}.
\end{proof}

 We shall now see that this result is a way to state Gelfand duality in
the real case.

 For this, we define first the upper real $\|a\|_0$ by:\[
\|a\|_0<r \text{ iff } 0\ll r-a \text{ and } 0\ll r+a.
                                                       \]
This defines a seminorm on $B$ which satisfies
$\|a^2\|_0 = \|a\|_0^2$; see~\cite{coquand05}.

 Each element $a$ defines a map of locales $\hat{a}:\MFN (A)\rightarrow
{\mathbb{R}}$
by taking $\hat{a}^{-1}(r,s)$ to be the open $D(a-r)\wedge D(s-a)$.
We define $\|\hat{a}\|$ as the upper real such that $\|\hat{a}\|<r$ iff
$1 = D(r-a)\wedge D(a+r)$.

\begin{proposition}\label{Gelfand1}
$\|\hat{a}\| = \|a\|_0$.
\end{proposition}
\begin{proof}
 By Proposition \ref{main}, $1 = D(r-a)\wedge D(a+r)$ is equivalent to
$0\ll a-r$ and $0\ll a+r$.
\end{proof}

\begin{corollary}\label{example}
$\|a\|^2_0 = \|a^2\|_0$.
\end{corollary}

\begin{proof}
This follows  from $\|\hat{a^2}\| = \|\hat{a}\|^2$ and Proposition
\ref{Gelfand1}.

Since Proposition \ref{main} is a combination of Lemma \ref{key} and
cut-elimination, we can also expect a direct proof from Lemma \ref{key}.
Here is such a direct argument. If $0\leqslant r$ and $0\ll r^2-a^2$ then
we have $0\ll uv$ where $u = r-a,~v = r+a$. Hence $0\ll u(u+v)$ and
$0\ll v(u+v)$. Since
$0\leqslant 2r = u+v$ we can apply Lemma \ref{key} and deduce $0\ll r+a$
and $0\ll r-a$.
\end{proof}

 To get Gelfand duality in the real case, we need to establish that $\|a\|_0$
and $\|a\|$ coincide. As usual the Stone-Weierstrass Theorem, which has a
constructive proof \cite{BMStone,coquand05}, then establishes the
surjectivity of
the map $a\longmapsto \hat{a}$.

\begin{lemma}\label{square}
$\|a^2 \| \leqslant  \|a^2\|_0$.
\end{lemma}

\begin{proof}
  Suppose that $\|a^2\|_0 < r$, then $r - a^2$ is a square, $b^2$. So
  \[\| a^2\|\leqslant \|a^2+b^2\| = r.\]
\end{proof}

\begin{theorem}\label{gelfand2}The Gelfand transform is norm-preserving:
 $\|a\|_0 = \|a\| = \|\hat{a}\|$.
\end{theorem}

\begin{proof}
We have $\|a\|_0\leqslant \|a\|$ since $r-a$ is a square if $r\geqslant
\|a\|$ by
Lemma \ref{lem:square}. On the other hand, we have
$\|a\|^2 = \|a^2\|\leqslant \|a^2\|_0 = \|a\|_0^2$ by
Corollary~\ref{example}
and Lemma~\ref{square}. Hence the result.
\end{proof}

\section{\label{gelfand}Constructive Gelfand duality}

 We now have all the pieces for constructive proof of Gelfand duality, also
in the complex case. Let $A$
be a commutative $C^*$-algebra and $B = A_{sa}$ its self-adjoint part. The
locale
$\MFN (A)$ defined in \cite{banaschewskimulvey00a} is isomorphic to the locale
$\MFN (B)$
defined above by
interpreting the element $a_1+ia_2\in (r_1+ir_2,s_1+is_2)$ in $\MFN (A)$ by the
element
$$D(a_1-r_1)\wedge D(s_1-a_1)\wedge D(a_2-r_2)\wedge D(s_2-a_2)$$
in $\MFN (B)$.

 Each element $b$ of $B$ defines a map of locales $\hat{b}:\MFN (A)\rightarrow
\mathbb{C}$
by taking $\hat{b}^{-1}(r,s)$ to be $b\in (r,s)$.

\begin{theorem}\label{gelfand3}
The Gelfand transform is norm-preserving: $\|b\| = \|\hat{b}\|$.
\end{theorem}

\begin{proof}
This follows from Theorem \ref{gelfand2}.
\end{proof}

\section{\label{applications}Some simple applications}

 We give some instances of simple properties of $C^*$-algebras that are
proved in \cite{banaschewskimulvey06} by using Barr's Theorem. All these
cases are direct
consequences of Proposition~\ref{sum:square} and do not depend on
Proposition~\ref{main}.

\begin{proposition}
  \label{oneminus}If $\|a\| \leqslant 1$, then $\|1 - a^{\ast} a\|
\leqslant 1$.
\end{proposition}

\begin{proof}
  Suppose that $\|a\| \leqslant 1$. Then $\|a^{\ast} a\| \leqslant 1$. Write
  $a = b + c i$, where $b, c$ are the real and the complex part. Then $a^{\ast}
a
  = b^2 + c^2$. Since $b^2 + c^2$ is a square it suffices to prove: If $\|d^2
  \| \leqslant 1$, then $\|1 - d^2 \| \leqslant 1$. Suppose that
$\|d^2 \|
  \leqslant 1$. Then $1 - d^2 = e^2$, so $1 = d^2 + e^2$ and hence $\|1 -
d^2
  \|=\|e^2 \| \leqslant 1$.
\end{proof}

\begin{proposition}
  The absolute value $( \sqrt{a^{\ast} a})$ exists.
\end{proposition}

\begin{proof}
  We can assume $\|a\|\leqslant 1$. Then $\|1 - a^{\ast} a\| \leqslant 1$. The
result now follows from
Lemma~\ref{lem:square}: the Taylor series $\sqrt{a^{\ast} a}$ converges.
\end{proof}

\begin{proposition}
Let $a$ be in $A$. Then $1 + a^{\ast} a$ is invertible.
\end{proposition}

\begin{proof}
  As in \cite{johnstone82}. Let $b^2 = a^{\ast} a$. Choose $n \geqslant 1 +
  b^2$. Define $c = (1 - \frac{1}{n}) - \frac{b^2}{n}$. By
  Proposition~\ref{oneminus}, $\|c\| \leqslant 1 - \frac{1}{n}$. It
follows
  that $(1 - c)^{- 1} = 1 + c + c^2 + \ldots$ exists and $n (1 - c)^{- 1}$ is
  the inverse of $1 + b^2$.
\end{proof}


\bibliographystyle{alpha}
\bibliography{gelfand.bib}

\end{document}